\documentclass[12pt]{article}

\def\edo{\end{document}	 }

\def\textgoth{\frak}

\usepackage{amsmath, amscd, amsfonts, amssymb}

\usepackage{appendix}                                 

\newtheorem{theorem}{Theorem}[section]

\newtheorem{corollary}[theorem]{Corollary}
\newtheorem{lemma}[theorem]{Lemma}
\newtheorem{definition}[theorem]{Definition}

\newtheorem{remark}[theorem]{Remark}

\usepackage[usenames,dvipsnames]{color}

\def\rrd{{\mathbb{R}^d}}
\def\call{{\mathcal{L}}}

\def\calr{{\mathcal{R}}}

\def\calf{{\mathcal{F}}}
\def\calm{{\mathcal{M}}}
\def\calo{{\mathcal{O}}}

\def\cals{{\mathcal{S}}}
\def\calb{{\mathcal{B}}}
\def\cald{{\mathcal{D}}}
\def\calx{{\mathcal{X}}}

\def\call{{\mathcal{L}}}

\def\calp{{\mathcal{P}}}

\def\pf{\n{\bf Proof.} }
\def\vsp{\vspace*{1,5mm}\\ }
\def\bk{\bigskip }
\def\mk{\medskip }
\def\sk{\smallskip }
\def\n{\noindent }
\def\dd{\displaystyle}

\def\barr{\begin{array}}
\def\earr{\end{array}}

\def\bit{\begin{itemize}}
\def\eit{\end{itemize}}

\def\FP{Fokker--Planck}

\def\1{^{-1}}
\def\one{\mbox{$1\!\!\,\rule{0,2mm}{3mm}\,$}}

\def\E{{\mathbb{E}}}
\def\rr{{\mathbb{R}}}
\def\nn{{\mathbb{N}}}
\def\9{{\infty}}
\def\lbb{{\lambda}}
\def\wt{\widetilde}

\def\vf{{\varphi}}

\def\ooo{{\Omega}}
\def\pp{{\partial}}
\def\vp{{\varepsilon}}

\def\ff{\forall }
\def\({\left(}
\def\){\right)}
\def\<{\left<}
\def\>{\right>}

\def\NS{Navier--Stokes}
\def\divv{{\rm div}}
\def\curl{{\rm curl}}
\def\D{{\Delta}}

\title{Uniqueness of distributional solutions to~the~2D vorticity   Navier--Stokes equation and its associated nonlinear Markov process} 
\author{Viorel Barbu\thanks{Al.I. Cuza University and Octav Mayer Institute of Mathematics of  Romanian Academy, Ia\c si, Romania.  Email: vbarbu41@gmail.com}\and Michael R\"ockner\thanks{Fakult\"at f\"ur Mathematik, Universit\"at Bielefeld,  D-33501 Bielefeld, Germany.  Email: roeckner@math.uni-bielefeld.de}\and Deng Zhang\thanks{School of Mathematical Sciences, Shanghai Jiao Tong University, Shanghai, China. Email: dzhang@sjtu.edu.cn}}
\date{}

\begin{document}
\maketitle\vspace*{-4mm}
\begin{abstract}
\n In this work we prove uniqueness of distributional solutions to~$2D$ Navier--Stokes equations in vorticity form $u_t-\nu\Delta u+\divv(K(u)u)=0$ on $(0,\9)\times\rr^2$ with  Radon measures as initial data, where $K$ is the Biot--Savart operator in 2-D. As a consequence, one gets the uniqueness of probabilistically weak solutions to the corresponding McKean--Vlasov stochastic differential equations. It is also proved that for initial conditions with density in $L^4$ these solutions are strong, so can be written as a functional of the Wiener process, and that pathwise uniqueness holds in the class of weak solutions, whose time marginal law densities are in $L^{\frac43}$ in space-time. In particular, one derives a  stochastic representation of the vorticity $u$ of the fluid flow in terms of a solution to the McKean--Vlasov SDE. Finally, it is proved that the family $\mathbb{P}_{s,\zeta},$ $s\ge0$, $\zeta=$probability measure on $\rrd$, of path laws of the solutions to the McKean--Vlasov SDE, started with $\zeta$ at $s$, form a nonlinear Markov process in the sense of McKean. \sk
\end{abstract}\vfill\newpage
\section{Introduction}\label{s1}
Consider here the 2-D incompressible \NS\ equation
\begin{equation}\label{e1.1}
\barr{ll}
y_t-\nu\D y+(y\cdot\nabla)y=\nabla p&\mbox{in }(0,\9)\times\rr^2,\vsp
\nabla\cdot y=0&\mbox{in }(0,\9)\times\rr^2,\vsp
y(0,x)=y_0(x),\ x\in\rr^2.
\earr\end{equation} 
\n Let $u=u(t,x)$ denote the vorticity of the velocity field $y=\{y_1,y_2\}$, that is,
$$u(t,x)=\curl\,y(t,x)=D_1y_2(t,x)-D_2y_1(t,x),\ (t,x)\in(0,\9)\times\rr^2,$$where $D_j=\frac\pp{\pp x_j},\ j=1,2,$ and the symbol  $\nabla,\divv$ refer to spatial derivatives. 

Equation \eqref{e1.1} can then be   rewritten as the {\it vorticity equation}
\begin{equation}\label{e1.2}
	\barr{ll}
	u_t-\nu\D u+\divv(yu)=0&\mbox{in }(0,\9)\times\rr^2,\\[1mm]
	u(0,x)=u_0(x)=\curl\,y_0(x),&x\in\rr^2.\earr\end{equation}
Here, the velocity field $y(t,x)$ is given by the Biot--Savart formula
\begin{equation}\label{e1.4}
y(t,x)=(\nabla^\bot E*u(t))(x),\ \ff(t,x)\in(0,\9)\times\rr^2,\end{equation}where $$E(x)=\frac1{2\pi}\,\ln|x|,\ x\in\rr^2,$$
hence
$$\nabla^\bot E(x)=\frac{(-x_2,x_1)}{2\pi|x|^2},\ x=(x_1,x_2)\in\rr^2\setminus\{0\},$$which is the Biot--Savart kernel. 
We set
\begin{equation}\label{e1.5}
K(z)=\nabla^\bot E*z,\ z\in L^p(\rr^2),\ p\in(1,2)\end{equation}
and note (see, e.g., \cite{11}, Lemma 2.2) that by the generalized Young inequality
\begin{equation}\label{e1.6a}
|K(z)|_{L^q(\rr^2)}\le C|z|_{L^p(\rr^2)},\ \ff z\in L^p(\rr^d), \end{equation}
and, if $y\in L^p(\rr^2,\rr^2)$ with ${\rm div}\,y=0$ and $u={\rm curl}\,y\in L^q$, that we have
\begin{equation}\label{e1.6}
y=K(u),\end{equation}	
where $p\in(1,2)$ and $\frac1q=\frac1p-\frac12.$  
Then, we may write \eqref{e1.2} as

\begin{equation}\label{e1.7}
\barr{ll}
u_t-\nu\D u+\divv(uK(u))=0&\mbox{in }(0,\9)\times\rr^2,\vsp
u(0,x)=u_0(x),&x\in\rr^2.\earr\end{equation} 
This is a special case of a so called {\it generalized mean-field \FP\ equation} with locally integrable singular kernel $K$ usually derived from Riesz potentials. Besides this kernel, other singular potentials of this form arise in chemotaxis mathematical models and in some classes of mean field equations (see    \cite{12a}--\cite{17a}, \cite{13a}).

There is an extensive literature on the well-posedness of the vorticity equation \eqref{e1.7} and, implicitly, on the Navier--Stokes equation \eqref{e1.1} in the spaces $L^p((0,\9)\times\rr^2)$ (see, e.g., \cite{6a}, \cite{11a}, \cite{11}, \cite{12}) and in $(BMO)\1(\rr^2)$, (see \cite{21a}).  

Our main objective here is the   relationship of \eqref{e1.7} with the McKean--Vlasov stochastic dif\-fe\-ren\-tial equation 
\begin{equation}\label{e1.8}
\barr{l}
dX(t)=K(u(t,\cdot))(X(t))dt+\sqrt{2\nu}\ dW(t),\ t\ge0,\vsp 
X(0)=X_0,\earr\end{equation}
on a probability space $(\ooo,\calf,\mathbb{P})$ with the normal filtration $(\calf_t)_{t\ge0}$ and 2-D $(\calf_t)$-Brownian motion $W(t)$, where
\begin{equation}\label{e1.9}
u(t,x)dx=\mathbb{P}\circ X(t)\1(dx);\ t>0,\  u_0(dx)=\mathbb{P}\circ X^{-1}_0(dx).\end{equation}SDE \eqref{e1.8} describes the microscopic dynamics of the vorticity flow $u=u(t,x)$ and, implicitly, of the Navier--Stokes velocity field $y(t)=K(u(t))$. Whereas weak existence of solutions to \eqref{e1.8} is a  consequence by the above mentioned existence results for \eqref{e1.7} (see, e.g., Theorem \ref{t2.1} below) and the general technique from Section 2.2 in \cite{2} (see Theorem \ref{t4.1} below), weak uniqueness results for \eqref{e1.8} appear to be less known.

Another open question is whether the path laws of the solutions to \eqref{e1.8} form a nonlinear Markov process in the sense of McKean (see \cite{24prim} and the recent paper \cite{27prim}). It turns out that to solve both problems, weak uniqueness for \eqref{e1.8} and the question whether the laws of its solutions lead to a Markov process, require a new uniqueness result for \eqref{e1.7}, namely uniqueness in the most general class of solutions for \eqref{e1.7}, the so-called distributional solutions (see \eqref{e2.2} and \eqref{e2.3} below). 

So, the first main result (which is purely analytic) of this paper is such a uniqueness result in the class of distributional solutions to \eqref{e1.7} with measure initial data, formulated as Theorem \ref{t2.2} (under the restrictions \eqref{e2.6}--\eqref{e2.7b}) and as Corollary \ref{c3.1}, where the latter is devoted to the case of probability measures as initial data. Furthermore, we also prove distributional uniqueness for the linearized equation corresponding to \eqref{e1.7} (see Theorem \ref{t4.2} and Corollary \ref{c4.3}, which is crucial for the subsequent probabilistic applications). As a consequence of both, we prove the existence of probabilistically weak solutions to \eqref{e1.8} (see Theorem \ref{t4.3}). For initial conditions which are pro\-ba\-bi\-li\-ty measures with densities in $L^4$ we prove that these solutions are in fact strong and that pathwise uniqueness holds in the class of all solutions whose path laws have time marginals densities in $L^{\frac43}$ in space-time (see Theo\-rem \ref{t4.5}). Our last main result, Theorem \ref{t4.6}, guarantees that the family of the path laws $\mathbb{P}_{s,\zeta}$ of the solutions to \eqref{e1.8}, started at time $s\ge0$ with pro\-ba\-bi\-lity measure $\zeta$ on $\rrd$, form a nonlinear Markov process in the sense of McKean \cite{24prim} (see Definition \ref{d4.7} and also \cite{27prim}). We would like to stress that all these results are heavily depending on our uniqueness result in Theorem \ref{t2.2}, and the fact that it gives uniqueness in the class of distributional solutions. Uniqueness in smaller classes as, e.g., mild solutions (see \cite{11} or the more general results in \cite{21prim}) does not suffice. We refer, e.g., to the proof of Theorem \ref{t4.3}, where this becomes obvious.

As a by-product of our existence result, Theorem \ref{t4.1}, for \eqref{e1.8}, we get a a probabilistic representation of the solutions to \eqref{e1.7} as time-marginal law densities of the nonlinear Markov process gotten from the paths laws of the solutions to \eqref{e1.8}. Thus, McKean's general   programme, already envisioned in \cite{24prim}, is completed in this paper for the $2D$ vorticity Navier--Stokes equation~\eqref{e1.7}. 

For the existence theory for nonlinear \FP\ equations with Ne\-mytski-type drift term and their implications to McKean--Vlasov SDEs, we refer to \cite{2}--\cite{4}. As regards the literature on the stochastic representation of solutions to \NS\ equations, the works \cite{8b}, \cite{11b}, \cite{13b}, \cite{12a} should be primarily cited. In particular, in \cite{8b} one gets the probabilistic representation of solutions to the vorticity equation \eqref{e1.7} as
$$u(t,x)=\E[u_0(X^{t,x}(t))],\ (t,x)\in(0,\9)\times\rr^d,$$where $X^{t,x}$ is the stochastic flow map defined by the equation
$$\barr{rcl}
dX^{t,x}(s)&=&-v(t-s,X^{t,x}(s))ds+\sqrt{2\nu}\,dW_s,\\
X^{t,x}(0)&=&x,\\
v(t,x)&=&-\dd\frac12\int^\9_0\frac1s\,\E[u(t,x+W(s))W^\bot_s]ds,\earr$$\newpage\n and $W(t)=(W^1(t),W^2(t)))$ is a $2{-}D$ Brownian motion and $W^\bot (t)=(W^2(t),W^1(t))$. 
This representation formula is extended later on in \cite{11b} to $3{-}D$ equations of the form \eqref{e1.1}. In the present paper, however,   we   use a different approach which takes advantage of the interpretation of the vor\-ti\-ci\-ty equation \eqref{e1.7} as a nonlinear \FP\ equation, which is associated with a McKean--Vlasov SDE by virtue of the superposition principle (see \cite{2} and see also \cite{13} for the case of usual SDEs).

Finally, let us comment on probabilistic approaches to \eqref{e1.8} and also \eqref{e1.7}. We start with pointing out the recent paper \cite{17a} which  contains a substantial discussion on the related literature, and we also refer to the re\-fe\-ren\-ces therein. Furthermore, we would like to mention  reference \cite{32} where also weak existence of solutions to \eqref{e1.8} for every probability measure on $\rr^2$ as initial data was proved, however, by a completely different method, not employing the nonlinear superposition principle. Furthermore, in \cite{17a} (see Theorem 6.3) it is proved that there exists a unique strong solution to \eqref{e1.8} with initial data in certain Besov spaces (which include $L^{1+\vp}$ data). So, this result is quite different from ours and again the methods are completely different from those in our present paper. In addition, uniqueness of smooth solutions to \eqref{e1.7} under additional restrictions on the behaviour at $t=0$ is proved \mbox{(see~\cite[Theorem 6.1]{17a}).} A further related paper should be mentioned, namely \cite{33}, in which a very nice theory for existence and uniqueness of soloutions to ordinary SDEs with singular drifts beyond the Ladyzhenskaya--Prodi--Serrin condition is developed, which applies to  \eqref{e1.8} after deter\-mi\-ning and fixing $u$. However, the general weak uniqueness results in \cite{17a} only give uniqueness in a class of martingale solutions obtained by a certain limiting procedure, and the Markov property is only proved for Lebesgue almost all times. So, our result in the present paper on \eqref{e1.8} giving rise to a  nonlinear Markov process (in  the sense of McKean, i.e., Theorem 4.6) is much more general.

\bk\n{\bf Notations.}  $L^p(\rr^2),\ 1\le p\le\9$ (denoted $L^p$) is the    space of all Lebesgue mea\-su\-rable and $p$-integrable functions on $\rr^2$, with the standard norm $|\cdot|_p$. $(\cdot,\cdot)_2$ denotes the inner product in $L^2$. 
By $L^p_{\rm loc}$ we denote the correspon\-ding local space.  For any open set $\calo\subset\rr^2$ let $W^{k,p}(\calo)$, $k\ge1$, denote\break the standard Sobolev space on  $\calo$ and by $W^{k,p}_{\rm loc}(\calo)$ the cor\-res\-pon\-ding local space. We set $W^{1,2}(\calo)=H^1(\calo)$,
$W^{2,2}(\calo)=H^2(\calo)$, $H^1_0(\calo)=\{u\in H^1(\calo),$ $u=0$ $\mbox{on }\pp\calo\}$, where $\pp\calo$ is the boundary of $\calo$. By $H\1=H\1(\rr^2)$ we denote the dual space of $H^1(\rrd)$. $C^\9_0(\calo)$ is the space of infinitely diffe\-ren\-tiable real-valued functions with compact support in $\calo$ and $\cald'(\calo)$ is the dual of $C^\9_0(\calo)$, that is, the space of Schwartz distributions on $\calo$. By~$C_b(\rr)$, we denote the space of continuous and bounded functions on $\rr^2$.  We shall denote by $\calm(\rr^2)$ the space of all finite Radon measures on $\rr^2$. Given a Banach space $\calx$ and $0<T\le\9$, we denote by $C([0,T];\calx)$ the space of all   continuous $\calx$-valued functions on $[0,T]$. For \mbox{$1\le p\le\9$,} we shall denote by $L^p(0,T;\calx)$ the space of $\calx$-valued, $L^p$-Bochner integrable functions on $(0,T)$. By  $C^\9_0([0,\9)\times\rr^2)$ we denote the space  $\{y\in C^\9([0,\9);\rr^2);$ $ y\mbox{ with compact support in }[0,\9)\}$.  For $1<p<\9$, let $L^{p,\9}(\rr^2)$ denote the Lorentz space of all measurable functions $y:\rr^2\to\rr^2$ such that 
$$|y|_{L^{p,\9}}:=\sup_{\lbb>0}\{\lbb^p\ {\rm meas}(x\in\rr^2;\ |y(x)|>\lbb)\}^{\frac1p}<\9.$$ 
Throughout this work, $\nabla={\rm grad}$ refers only to the spatial derivatives, i.e., in the $x$-variables and 
$$\nabla\cdot y=\divv\,y,\ \ff y\in(L^p(\rr^2))^2,\ 1\le p\le\9.$$If $X_1,X_2$ are two Banach spaces, we shall denote by  $L(X_1,X_2)$ the space of linear continuous operators from $X_1$ to $X_2$. 
By $\calp$ we denote the set of all probability measures on $\rr^2$ and set
\begin{equation}\label{e1.10}
	\calp^a=\left\{y\in L^1(\rr^2);y\ge0,\mbox{ a.e. in }\rr^2;\int_{\rr^2}y(x)dx=1\right\}.
	\end{equation}

\section{The existence and uniqueness\\ for the vorticity equation \eqref{e1.7}}\label{s2}
\setcounter{equation}{0}

A function $u\in L^{r_1}_{\rm loc}(0,\9;L^{r_2}),\ r_1,r_2\ge 1$, is called a {\it mild solution to \eqref{e1.7}} if it is a solution to the integral equation 
		
\begin{equation}\label{e2.1}
u(t)=e^{\nu t\Delta}u_0-\divv\int^t_0 e^{\nu(t-s)\D}(K(u(s))u(s))dx,\ t>0,\end{equation}
where $e^{t\D}$ is the heat semigroup in $\rr^2$, which is well defined on all $L^p$,\break $1\le p\le\9$. (If $u_0\in\calm(\rr^2)$, then $\|e^{t\D}u_0\|_{L^p}\le C t^{-1+\frac1p}\|u_0\|_\calm,$ $\ff t>0$, for all $p>1.$) 

 This definition extends to mild solutions $y$ to \eqref{e1.1} via the Biot--Savart formula \eqref{e1.6}. 

Given $u_0\in\calm(\rr^2)$, a function $u\in L^1_{\rm loc}(0,\9;L^p)$ for some $p\in(1,2)$ is said to be a distributional solution to \eqref{e1.7} if
 \begin{eqnarray}
 &K(u)u\in L^1_{\rm loc}((0,\9)\times\rr^2),\label{e2.2}\\[2mm]
 &\barr{r}
 \dd\int^\9_0\!\!\int_{\rr^2}u(t,x)(\vf_t(t,x)+\nu\D\vf(t,x)+K(u(t,x))\cdot\nabla\vf(t,x))dtdx\\
 +\dd\int_{\rr^2}\vf(0,x)u_0(dx)=0,\ \ff\vf\in C^\9_0([0,\9)\times\rr^2).\earr\label{e2.3}
\end{eqnarray}
 
The next existence theorem is due to Y. Giga, T. Miyakawa \& H. Osada (see \cite{11}, Theorem 4.2).

\begin{theorem}\label{t2.1} Assume that $u_0=\curl\,y_0\in\calm(\rr^2)$, where $y_0\in L^{2,\9}(\rr^2),$ $\nabla\cdot y_0=0$ on $\rr^2$. Then, equation \eqref{e1.7} has a solution $u:[0,\9)\to\calm(\rr^2)$, which is bounded and continuous in the weak topology. Moreover, one has the estimate
\begin{equation}\label{e2.4}
|u(t)|_r\le C_r\,t^{-1+\frac1r},\ \ \ff\,t>0,\ 1<r<\9,	
\end{equation}
	and $u(t)=\curl\,y(t)$, where $y$ is a mild solution to the Navier--Stokes equation \eqref{e1.1} and
\begin{equation}\label{e2.5}
	|y(t)|_r\le C_r\,t^{\frac1r-\frac12},\ \ \ff\,t>0, 2<r\le\9,
\end{equation}
	\begin{equation}\label{e2.5a}
\sup\{|D_x^kD^j_tu(t)|_\9;\ t\in[\vp,T]\}\le C^T_{\vp,k,j}, 
\end{equation}
for all $0<\vp<T<\9$ and all $k,j=0,1,...$ 
\end{theorem} 

\begin{remark}\label{r2.1prim}\rm According to Theorem 1.2 in \cite{18a}, Theorem \ref{t2.1} extends to all initial conditions $u_0\in\calm(\rr^2)$.\end{remark}

We set
\begin{equation}\label{e2.6'}
\barr{r}
k(x)=(k^1(x),\ k^2(x)):=\nabla^\bot E(x)=
\dd\frac{(-x_2,x_1)}{2\pi|x|^2},\vsp 
x\in(x_1,x_2)\in\rr^2\setminus\{0\},\earr
\end{equation}and
$$K^i(u):=k^i*u,\ u\in L^p,\ p\in(1,2),\ i=1,2.$$
We have
$$\nabla k(x)=(\pp_jk^i)_{1\le i,j\le2}=\frac1{2\pi|x|^2}\(
\barr{ccc}
\dd\frac{2x_1x_2}{|x|^2}&&-\(1-\dd\frac{2x^2_2}{|x|^2}\)\vsp 
\(1-\dd\frac{2x^2_1}{|x|^2}\)&&-\dd\frac{2x_1x_2}{|x|^2}\earr\),$$and hence, all four components $\pp_jk^i$ of $\nabla k(x)$ define kernels of degree $(-2)$ satisfying all assumptions of Theorem 1 in \cite{10}, which implies that, for all $p\in(1,\9)$; $i,j\in\{1,2\},$ 
\begin{equation}\label{2.6prim}
\pp_jK^i(u)(x):=\lim_{\vp\to0+}
\int_{|x-y|\ge\vp}\pp_j k^i(x-y)u(y)dy,\ x\in\rr^2,\ u\in L^p,\end{equation}
defines bounded, linear operators $\pp_i K^j:L^p\to L^p.$ The limit in \eqref{2.6prim} is meant in $L^p$ as well as a.e.  Moreover, the $p$-norm of the right hand side of  \eqref{2.6prim} with $\sup\limits_{\vp\in(0,1]}$ replacing $\lim\limits_{\vp\to0+}$ is up to a constant (only depending on $p$) dominated by $|u|_p$. On the other hand, it is elementary to check and well known that for the distributional derivatives $\pp_j(K^i(u))$ of $K^i(u)$, $i=1,2,$ $u\in L^p$, $p\in(1,2),$ 
\begin{equation}
\label{e2.7prim}
\pp_j(K^i(u))=\pp_jK^i(u)+\frac12\,{\rm sign}(i-j)u
\end{equation} for some numerical constant  $c>0.$ In particular,  we have
\begin{equation}\label{2.6tert}
{\rm div}(K(u))=0\mbox{ and }{\rm curl}(K(u))=u,\  \ff u\in L^p,\ p\in(1,2), \end{equation}where div and curl are taken in the sense of Schwartz distributions. Furthermore, together with the fact that the operators in  \eqref{2.6prim} are bounded on every $L^p$, $p\in(1,\9)$, 
\eqref{e2.7prim} implies that
\begin{equation}\label{2.6secund}
|\nabla K(u)|_p\le c_p|u|_p,\ \ff u\in L^p\cap\(\bigcup\limits_{q\in(1,2)}L^q\),\ p\in(1,\9). \end{equation}


Theorem \ref{t2.1} can be complemented  as follows.

\begin{theorem} \label{t2.2bis} Let $u$ be the solution of \eqref{e1.7} from Theorem {\rm\ref{t2.1}} and let $T>0$. Then $u$ is a distributional solution  to equation \eqref{e1.7} and
\begin{equation}\label{e2.9prim}
{\rm(i)}\hspace*{20mm} \int^T_0\int_{\rr^d}|u(t,x)|\,|K(u(t,\cdot))(x)|dx\,dt<\9.\hspace*{20mm}		\end{equation}  

\begin{itemize}\item[\rm(ii)] If $u_0\in\calp$, then
	\begin{equation}
	\label{2.7}
	u(t)\in\calp^a,\ \ff t>0.
	\end{equation}
\item[\rm(iii)] Assume that $u_0\in L^q$ for some $q\in(1,2)$. Then, for all $p\in[1,\9]$ we have $u\in L^\9(0,T;L^p)$ and 
\begin{equation}\label{2.12prim}
	|u|_{L^\9(0,T;L^p)}\le c|u_0|_p,\end{equation}
where $c$ only depends on $p$. \end{itemize}

Therefore, $u\in L^4(0,T;L^{\frac43}\cap L^4)\cap L^\9(0,T;H\1)$ for $u_0\in L^{\frac43}\cap L^4$.
\end{theorem}

\pf (i): Let $p\in(1,2)$, $q:=\frac{2p}{2-p}$ and $q':=\frac q{q-1}$. Then, by \eqref{e1.6a},
$$\int^T_0\int_{\rr^d}|u(t,x)|\,|K(u(t,x))|dx\,dt
\le C\int^T_0|u(t)|_{q'}|u(t)|_pdt.$$
But by \eqref{e2.4} the last integral is finite. In particular, \eqref{e2.2}  holds, hence clearly $u$ is a distributional solution to \eqref{e1.7}.\mk  

\n(ii): As seen in the proof of Theorem 4.2 in \cite{11}, the solution $u$ to \eqref{e1.7} given by Theorem \ref{t2.1} can be obtained by 
\begin{equation}\label{2.8}
u=\lim_{u\to\9} u_n\mbox{ uniformly on compact sets of $(0,\9)\times\rr^2$}\end{equation}where $u_n$ are the unique global smooth solutions to \eqref{e1.7} with $u_n(\cdot,0)=u^n_0$ and $\{u^n_0\}$ is a smooth approximation of the initial data $u_0\in\calm(\rr^2)$ in the narrow topology. Moreover, $u_n$ is expressed as (see, e.g., (2.1) in \cite{11})
\begin{equation}\label{2.8a}
	u_n(t,x)=\int_{\rr^2}\Gamma_{K(u_n)}(t,x;0,\xi)u^n_0(\xi)d\xi,\end{equation}
where $\Gamma_v\equiv\Gamma_v(t,x;s,\xi)$  is the fundamental solution to the linear parabolic operator
$$L_v(u)=u_t-\nu\Delta u+(u\cdot\nabla)v,\ (t,x)\in(0,\9)\times\rr^2.$$ 
 We have
\begin{equation}\label{2.14a}
\barr{c}
\Gamma_v\ge0;\ \dd\int_{\rr^2}\Gamma_v(t,x;s,\xi)d\xi=1,\ 0\le s\le t<\9,\ x\in\rr^2,\vsp 
\dd\lim_{t\downarrow s}\int_{\rr^2}\Gamma_v(t,x;s,\xi)f(\xi)d\xi=f(x),\ \ff f\in C_b(\rr^2).\earr\end{equation}
If $u_0\in\calp$, then the sequence $\{u_0^n\}$ can be chosen in such a way that $u^n_0\ge0$ and $u^n_0\to u_0$ in $\calm(\rr^2)$ narrowly as $n\to\9$. Then, by  \eqref{2.8}--\eqref{2.14a} it follows that $u\ge0$. 
Furthermore, by Theorem 1.2 in \cite{18a} we know that  $$\int_{\rr^2}u(t,x)dx=\int_{\rr^2}u_0(x)dx=1,\ \ff t\in(0,T),$$as claimed.\mk

\n(iii) Since $q\in(1,2)$, by \eqref{2.6tert} we know that $u={\rm curl}(K(u))$. Hence, we may apply Theorem 4.3 in \cite{11} to obtain the following representation
\begin{equation}\label{2.16'}
	u(x,t)=\int_{\rr^2}\Gamma(t,x;0,\xi) u_0(x)d\xi,\ t\ge0,\ x\in\rr^d,\end{equation}where $\Gamma(t,x;s,\xi)$, $x,\xi\in\rr^2,\ t>s\ge0$, is a positive continuous function, which satisfies
$$\barr{c}
\dd\int_{\rr^2}\Gamma(t,x;s,\xi)d\xi=\int_{\rr^2}\Gamma(t,x;s,\xi)dx=1,\vsp 
\Gamma(t,x;s,\xi)\simeq\dd\frac1{(t-s)}\ e^{-\frac{C|x-\xi|^2}{(t-s)}},\ t>s\ge0.\earr$$
	By the Young inequality $\ff 1\le p\le\9$, if $u_0\in L^p$, this yields
	$$\barr{ll}
	|u(t)|_{p}\!\!\!
	&=\left|\dd\int\Gamma(t,x;0,\xi)u_0(\xi)d\xi\right|_{p}\vsp 
&\le c \left|\dd\int\frac1{t-s}\ 
e^{-\frac{C|x-\xi|^2}{t-s}}\ u_0(\xi)d\xi\right|_{p}\vsp
&\le c \left|\dd\int\frac1t\ e^{-\frac{C|x|^2}t}\ dx\right||u_0|_p\vsp 
&\le c|u_0|_p.\earr$$Thus, $u\in L^\9(0,T;L^p)$, and $|u|_{L^\9(0,T;L^p)}\le c|u_0|_p$, as claimed.\hfill$\Box$\bk
 
Theorem \ref{t2.1}  is completed in \cite{11} by a uniqueness theorem for \eqref{e1.1} and implicitly for \eqref{e1.7}, in the class of mild solutions with sufficiently small atomic  part $(u_0)_{pp}$ of the Radon measure $u_0$. Such a uniqueness result was extended in \cite{18a} to $u_0\in \calm(\rr^2)$ in the class of mild solutions $u\in C((0,T];L^1\cap L^\9)\cap L^\9(0,T;L^1)$.

We shall prove here for the purposes of the McKean--Vlasov equation \eqref{e1.8} a sharper uniqueness result, namely within the much larger class of all  dis\-tri\-bu\-tional solutions to \eqref{e1.7}, in the sense of \eqref{e2.2}, \eqref{e2.3}, which belong to the class $\{u\!\in\! L^4(0,T;$ $L^4\cap L^{\frac43});\ \ff T>0\}$. More precisely, we prove  
 \begin{theorem}\label{t2.2} Let $u_1,u_2$ be two distributional solutions to \eqref{e1.7} such that  
 \begin{eqnarray}
&u_1,u_2\in L^4(0,T;L^4\cap L^{\frac43}),\ \ff T>0,\label{e2.6}\\[1mm]
&u_1-u_2\in L^\9(0,T;H\1),\label{e2.7}\\[1mm]
&u_1(t)-u_2(t)\in L^1\mbox{ for a.e. }t\in(0,T),\label{e2.7b}\\[1mm] 
 &\dd\lim_{t\downarrow0}\ {\rm ess}\sup_{\hspace*{-5mm}s\in(0,t)}\int_{\rr^2}(u_1(s,x)-u_2(s,x))\vf(x)dx=0,\ \ff\vf\in C^\9_0(\rr^2).\label{e2.8}
 		\end{eqnarray}Then, $u_1\equiv u_2$.
 \end{theorem}
	
We note that the uniqueness class considered here is larger than that covered by \cite{18a} and the method of proof is different. In fact, in the case of nonlinear \FP\ equations with Nemytski-type drift term, a result of this type was proved in \cite{5} in the class of $L^2((0,T)\times\rr^2)$ distributional solutions $u_1,u_2$ such that $u_1-u_2\in L^\9(0,T;H\1),$ $ \ff T>0$, but there is not a large overlap. Though the idea of the proof is borrowed from \cite{5}, the argument used here requires sharp estimates specific to the drift term ${\rm div}(uK(u))$.   It should be also mentioned that the uniqueness condition \eqref{e2.8} does not exclude the class of solutions with measure initial data $u_0\in\calm(\rr^2)$. In fact, such a condition agrees with Theo\-rem \ref{t2.1} which provides a distributional solution $u:[0,\9)\to\calm(\rr^2)$, which is weakly continuous.  We also note that, as in the proof of \eqref{e2.9prim}, by \eqref{e2.6} it follows that $u_i K(u_i)\in L^1_{\rm loc}((0,\9)\times\rrd),$ $i=1,2,$ which is a condition required by \eqref{e2.2}. 	

\section{Proof of Theorem \ref{t2.2}}\label{s3}
\setcounter{equation}{0}

We set $z=u_1-u_2.$ Then, we have by \eqref{e1.7}
\begin{equation}\label{e3.1}
z_t-\nu\D z+\divv(K(z)u_1-zK(u_2))=0\mbox{ in }\cald'((0,T)\times\rr^2),\end{equation}where $K$ is the Biot--Savart operator  \eqref{e1.5}. 


We also recall that  in 2-D we have
\begin{equation}\label{e3.4a}
	|w|^2_4\le2|w|_2|\nabla w|_2,\ \ff w\in H^1.\end{equation}It follows by \eqref{e1.6a} that, for $u_1,u_2\in L^4\cap L^{\frac43}$, 
\begin{equation}\label{e3.4}
|u_1K(z)|_2\le|K(z)|_4|u_1|_4\le 
 C|z|_{\frac43}|u_1|_4.\end{equation}Similarly, we get the estimate
\begin{equation}\label{e3.5}
|zK(u_2)|_2\le C|z|_4|u_2|_{\frac43},\end{equation}for some constant $C$ independent of $u_1,u_2$. Since for $u\in L^{\frac43}\cap L^4$, 
by interpolating between $L^{\frac43}$ and $L^4$,  it follows also that $|u|_2\le|u|^{\frac12}_{\frac43}|u|^{\frac12}_4$.  Therefore, since $u_1,u_2\in L^2(0,T;L^{\frac34}\cap L^4),$
\begin{equation}\label{e3.6}
z\in L^2(0,T;L^2)\cap L^\9(0,T;H\1).\end{equation} 
Consider the operator  $\Phi_\vp:L^2\to L^2,$  
\begin{equation}
\label{e3.6a}
\Phi_\vp(y)=(\vp I-\D)\1y,\ \ff y\in L^2,\ \vp>0.
\end{equation} 
and we note that 
\begin{equation}
\label{e3.9a}
\Phi_\vp\in L(L^2,H^2)\cap L(H\1,H^1)\cap L(L^2,L^2).\end{equation}
Then, applying $\Phi_\vp$ in \eqref{e3.1}, we get
\begin{equation}\label{e3.7}
\barr{r}
(\Phi_\vp(z(t)))_t{-}\nu\D\Phi_\vp(z(t)){+}\Phi_\vp(\divv(K(z(t))u_1(t){-}z(t)K(u_2(t))))=0\vsp\mbox{ in }\cald'((0,T)\times\rr^2).\earr\end{equation}Taking into account that, by \eqref{e3.4}, \eqref{e3.5}, and \eqref{e2.6},
$$K(z)u_1-zK(u_2)\in L^2(0,T;L^2),$$it follows by  \eqref{e3.9a}, \eqref{e3.7} and \eqref{e3.6} that
\begin{equation}\label{e3.8}
\frac d{dt}\ \Phi_\vp(z)\in L^2(0,T;L^2)\end{equation}and since, by \eqref{e3.6} and \eqref{e3.9a},  $\Phi_\vp(z)\in L^2(0,T;H^2)$, we infer that $\Phi_\vp(z)\in C([0,T];H^1)$. In particular, this implies that there is

\begin{equation}\label{e3.9}
\lim_{t\to0}\Phi_\vp(z(t))=f_\vp\mbox{ in }H^1.\end{equation}Now, we set
\begin{equation}\label{e3.10}
h_\vp(t)=(\Phi_\vp(z(t)),z(t))_2,\ t\in(0,T),\end{equation}and
\begin{equation}\label{e3.11}
K_\vp(y)=\nabla^\bot\Phi_\vp(y),\ \ff y\in L^2,\ \vp>0.\end{equation}We note that $K_\vp(z)=\nabla^\bot\Phi_\vp(z)\in L^2(0,T;H^1)\cap C([0,T];L^2)$. Taking into account that, by \eqref{e3.6a},
\begin{equation}\label{e3.11a}
\vp\Phi_\vp(z)-\D\Phi_\vp(z)=z\mbox{ on }\rr^2, \end{equation}it follows by \eqref{e3.10}, \eqref{e3.11} that
\begin{equation}\label{e3.12}
h_\vp(t)=|K_\vp(z(t))|^2_2+\vp|\Phi_\vp(z(t))|^2_2,\  t\in(0,T).\end{equation}Since $\Phi_\vp(z)\in C([0,T];H^1)$, by  \eqref{e3.8} we see that $h_\vp$ is absolutely continuous on $[0,T]$ and, by \eqref{e3.7} and \eqref{e3.11a} we have, for $\vp\in(0,1)$,
\begin{equation}\label{e3.13}
\hspace*{-4mm}\barr{ll}
h'_\vp(t)\!\!\!
&=2\(\dd\frac d{dt}\,\Phi_\vp(z(t)),z(t)\)_2\vsp
&=2(\nu\D\Phi_\vp(z(t)){-}\Phi_\vp(\divv(K(z(t))u_1(t){-}z(t)K(u_2(t)))),z(t))_2\vsp
&=-2\nu|z(t)|^2_2+2\vp\nu(\Phi_\vp(z(t)),z(t))_2\vsp
&\ \ \ 
+2(K(z(t))u_1(t)-z(t)K(u_2(t)),\nabla\Phi_\vp(z(t)))_2\vsp 
&\le-2\nu|z(t)|^2_2+2\vp\nu h_\vp(t)
+2(K(z(t)u_1(t),\nabla\Phi_\vp(z(t))))_2\vsp
&\ \ \ +2|z(t)K(u_2(t))\cdot\nabla\Phi_\vp(z(t))|_1,   \mbox{ a.e. }t\in(0,T).\qquad\, 
\earr\hspace*{-10mm}
\end{equation}On the other hand, we have 
\begin{equation}\label{e3.14}
h_\vp(0+)=\lim_{t\downarrow0}h_\vp(t)=0.\end{equation} Indeed, we have, for a.e. $t>0$,
$$\barr{r}
0\le h_\vp(t)\le |\Phi_\vp(z(t))-f_\vp|_{H^1}\,|z|_{L^\9(0,T;H\1)}+(f_\vp-\vf,z(t))_2+(\vf,z(t))_2,\vsp \ff\vf\in C^\9_0(\rr^2),\earr$$and so, by \eqref{e3.9} and \eqref{e2.8} 
we have
$$h_\vp(0+)=\dd\lim_{t\downarrow0}\ {\rm ess}\sup_{\hspace*{-5mm}s\in(0,t)}h_\vp(s)
\le |f_\vp-\vf|_{H^1}|z|_{L^\9(0,T;H\1)},\ \ff\vf\in C^\9_0(\rr^2).$$Since $|f_\vp-\vf|_{H^1}$ can be chosen sufficiently small for a suitable $\vf\in C^\9_0(\rr^2)$, by \eqref{e2.7} we get \eqref{e3.14}, as desired.\newpage

\n Now, taking into account that $|\nabla\Phi_\vp(z(t))|^2_2=|K_\vp(z(t))|^2_2$, we see by  \eqref{e3.13} that
\begin{equation}\label{e3.15}
	\barr{r}
h'_\vp(t){+}2\nu|z(t)|^2_2
\le 2\nu\vp h_\vp(t)+2|K_\vp(z(t))|_2\,|K(z(t))u_1(t)|_2\vsp 
 +2|z(t)K(u_2(t))\cdot\nabla\Phi_\vp(z(t))|_1,
 \mbox{ a.e. }t\in(0,T).\earr
\end{equation}Note that 
\begin{equation}\label{e3.15a}
\barr{r}
|K_\vp(z(t))|_2|K(z(t))u_1(t)|_2 
 \le|K(z(t))|_4|u_1(t)|_{4}|K_\vp(z(t))|_2, \vsp
  \le |K(z(t))|^2_4+|u_1(t)|^2_{4}|K_\vp(z(t))|^2_2,\ \mbox{ a.e. }t\in(0,T).\earr\end{equation}We also have, by \eqref{e1.6a}  and \eqref{e3.4a},     
 \begin{equation}\label{e3.15aa}
 \barr{l}
 \!\!\!|z(t)K(u_2(t))\cdot \nabla\Phi_\vp(z(t))|_1\le
 |z(t)|_2|K(u_2(t))|_4|K_\vp(z(t))|_4\vsp
 \quad\le C|z(t)|_2
 |u_2(t)|_{\frac43}|K_\vp(z(t))|^{\frac12}_2|\nabla K_\vp(z(t))|^{\frac12}_2\vsp
 \quad\le C|z(t)|_2|u_2(t)|_{\frac43}||K_\vp(z(t))|^{\frac12}_2|z(t)|^{\frac12}_2\vsp
 \quad\le C|z(t)|^{\frac32}_2
 |u_2(t)|_{\frac43}||K_\vp(z(t))|^{\frac12}_2\vsp
 \quad\le \dd\frac\nu4 |z(t)|^2_2+C_\nu
 |u_2(t)|^4_{\frac43}|K_\vp(z(t))|^2_2,\mbox{ a.e. }t\in(0,T).\earr\end{equation}
Here, we have used the inequality  
\begin{equation}\label{e3.23b}
|\nabla K_\vp(z)|_2\le C|z|_2,\ \ff z\in L^2\cap L^{\frac43},\end{equation}
which follows, since by Lemma \ref{A.1} in the Appendix we know that
$$K_\vp(z)=-K(z)+\vp K(\Phi_\vp(z)),$$and, therefore, by  \eqref{2.6secund},
$$\barr{ll}
|\nabla K_\vp(z)|_2\!\!\!
&\le |\nabla K(z)|_2+\vp|\nabla K(\Phi_\vp(z))|_2\vsp &\le C(|z|_2+\vp|\Phi_\vp(z)|_2)\vsp 
&\le 2C|z|_2.\earr$$
(Here and everywhere in the following we have denoted by $C$ several positive constants independent of $\vp$ and $u_1,u_2$.) 
 
Then, substituting   \eqref{e3.15a}--\eqref{e3.15aa} into \eqref{e3.15} and recalling that, by \eqref{e3.12}, $|K_\vp(z)|^2_2\le h_\vp$,  yields

\begin{eqnarray}
\!\!\!\!\!\!\!\!\!\!\!\!h'_\vp(t)&\!\!\!\!+\!\!\!\!&\nu|z(t)|^2_2\label{e3.16} \\
&\!\!\!\!\le\!\!\!\!& C(\vp h_\vp(t){+}|u_1(t)|^2_4|K_\vp(z(t))|^2_2 
{+}|K(z(t))|^2_4
 {+}|u_2(t)|^4_{\frac43}|K_\vp(z(t))|^2_2)\nonumber\\[1mm]
 &\!\!\!\! \le\!\!\!\!& C\(\(\vp{+}|u_1(t)|^2_4{+}|u_2(t)|^4_{\frac43}\)h_\vp(t)\){+}C|K(z(t))|^2_4,  
 \mbox{ a.e. }t\in(0,T),\nonumber 
\end{eqnarray}
we get by \eqref{e3.16} that
\begin{equation}\label{e3.23bb}
\barr{l}
\!\!\!\!\!\!\!\!\!\!\!\!\dd\frac d{dt}\!
\(\!h_\vp(t)\exp\!\(\!-C
\(\vp t+\!\!\int^t_0\!\!\(|u_1(s)|^2_4{+}
|u_2(s)|^{4}_{\frac43}\)ds\!\)\!\)\)\vsp
 \le\! C|K(z(t))|^2_4\exp\!
\(\!-C\!\(\vp t{+}\!\dd\int^t_0\!\!
\(|u_1(s)|^2_4+|u_2(s)|^{4}_{\frac43}\)ds\!\)\!\)\!,\vsp
\hfill \mbox{a.e. }t\in(0,T).\earr\end{equation}Since, by \eqref{e1.6a},
$$|K(z(t))|^2_4\le C(|u_1(t)|^2_{\frac43}+|u_2(t)|^2_{\frac43}), \mbox{ a.e. }t\in(0,T),$$and so $|K(z)|^2_4\in L^1(0,T)$, we see by \eqref{e3.14} and \eqref{e3.23bb} that
\begin{equation}\label{e3.26a}
\barr{ll}
0\le h_\vp(t)\vsp
\le C\!\!\dd\int^t_0\!\!\!|K(z(s))|^2_4
\exp\!\(\!C\!
\(\vp(t{-}s){+}\!\!\int^t_s\!\!
\(|u_1(\tau)|^2_4{+}|u_2(\tau)|^{4}_{\frac43}\)
d\tau\!\!\)\!\)ds,\vsp 
\hfill \ff\vp>0,\ t\in[0,T].\earr
\end{equation}Taking into account that, by \eqref{e2.6},
$$|u_1|^2_4+|u_2|^{4}_{\frac43}\in L^1(0,T),$$
by \eqref{e3.12} and \eqref{e3.26a} we have
\begin{equation}\label{e3.23a}
\sup_{\vp\in(0,1]}\|K_\vp(z)\|_{C([0,T];L^2)}
=\sup_{\vp\in(0,1]}\|\nabla^\bot(\vp I-\Delta)\1z\|_{C([0,T];L^2)}<\9.
\end{equation}
We set $\theta_\vp(t)=\Phi_\vp(z(t))$ and note that
\begin{equation}
\label{3.27a}
\vp\theta_\vp(t)-\D\theta_\vp(t)=z(t),\mbox{ a.e. }t\in(0,T)
\end{equation}
and, by \eqref{e3.12}, \eqref{e3.26a}, it follows that, for all $\vp\in(0,1),\ t\in[0,T]$,

\begin{equation}\label{3.28}
\vp|\theta_\vp(t)|^2_2+|\nabla\theta_\vp(t)|^2_2=h_\vp(t)
\le\sup_{^{\vp\in(0,1)}_{t\in[0,T]}}h_\vp(t)\le  C<\9.\end{equation}
We introduce the space 
$$\mathcal{G}:=\{u\in L^2_{\rm loc}(\rr^2):|\nabla u|\in L^2(\rr^2)\}$$equipped with the inner product $(\nabla\cdot,\nabla\cdot)_2$ (see p.~11 in \cite{17*}).  Then, by \cite{16*}, we have
\begin{itemize}
\item[$(\mathcal{G}.1)$] $\mathcal{G}=\left\{T\in\cald'(\rr^2):\frac{\pp T}{\pp x_i}\in L^2(\rr^2),\ 1\le i\le2\right\}.$ 
\item[$(\mathcal{G}.2)$] The quotient space $$\dot{\mathcal{G}}:=\mathcal{G}/\{constants\}$$
is a Hilbert space. Furthermore, for every Cauchy  sequence $u_n\in\mathcal{G},$ $n\in\nn$, there exist $u\in\mathcal{G}$ and $c_n\in\rr$ such that $\lim\limits_{n\to\9}u_n=u$ in $\mathcal{G}$
and $\lim\limits_{n\to\9}(u_n+c_n)=u$ in $L^2_{\rm loc}$.  
\end{itemize}
By \eqref{3.28} there exist subsequences  $\vp_k\in(0,1],$ $k\in\nn
$, and $\ell_n\in\nn$, $n\in\nn$, such that $\vp_k\to0,$ as $k\to\9$ and for $v_n:=\frac1{\ell_n}\sum\limits^{\ell_n}_{k=1}\theta_{\vp_k},\ n\in\nn,$ we have
\begin{equation}\label{e3.26'}
\nabla v_n\to F\mbox{ in } L^2((0,T)\times\rr^2),\end{equation}
and, for some Lebesgue zero set $N\subseteq(0,T)$,  
\begin{equation}\label{e3.27prim}
\nabla v_n(t)\to F(t)\mbox{ in }L^2(\rr^2), \mbox{ for every $t\in(0,T)\setminus N$.}\end{equation}Below we fix $t\in(0,T)\setminus N$. By $(\mathcal{G}.2)$ we know that there exist $\theta(t)\in\mathcal{G}$ and $c_n(t)\in\rr$ such that
\begin{equation}\label{3.30}
\lim_{n\to\9}\nabla v_n(t)=\nabla\theta(t)\mbox{ in }L^2
\end{equation}and
\begin{equation}\label{3.31}
\lim_{n\to\9}(v_n(t)+c_n(t))=\theta(t)\mbox{ in }L^2_{\rm loc}. 
\end{equation}
Furthermore, by \eqref{3.27a} we have, for every $n\in\nn$,
\begin{equation}\label{3.32}
\frac1{\ell_n}\sum^{\ell_n}_{k=1}\vp_k\theta_{\vp_k}(t)-\Delta v_n(t)=z(t).  
\end{equation}Hence, taking the limit $n\to\9$ in $\cald'(\rr^2)$, by  \eqref{3.28} we conclude that, for all $t\in(0,T)\setminus N$,
\begin{equation}\label{e3.30prim}
-\Delta\theta(t)=z(t),\  \mbox{ in }\cald'(\rr^2).\end{equation} 
 Moreover, by \eqref{3.30}--\eqref{3.31} it follows that $\theta(t)\subset W^{1,2}_{\rm loc}(\rr^2),$ for all \mbox{$t\in(0,T)\setminus N$,} and by \eqref{e3.26'}, \eqref{e3.27prim} that $\nabla\theta\in L^2((0,T)\times \rr^2).$ Furthermore, by \eqref{3.28}
$$|\nabla\theta(t)|_2\le C<\9,\ \mbox{ a.e.}t\in(0,T).$$This yields
$$\barr{l}
\dd\lim_{n\to\9}\int_{[1\le|x|\le2]}
|\nabla\theta(t,nx)|dx=\dd\lim_{n\to\9}\frac1{n^2}\int_{[n\le|y|\le2n]}|\nabla\theta(t,y)|dy\vsp
\qquad\dd\le\sqrt{3\pi}\lim_{n\to\9}\frac1n
\(\int|\nabla\theta(t,y)|^2dy\)^{\!\!\frac12}=0,\mbox{ a.e. }t\in(0,T).\earr$$
We recall also that $z(t)\in L^1$, a.e. $t\in(0,T)$ by \eqref{e2.7b}. Then, by Lemma A.11 in \cite{6b}, we have
\begin{equation}\label{e3.30secund}
\nabla\theta(t)=-\nabla E*z(t),\mbox{ a.e. }t\in(0,T),\end{equation}and this yields
\begin{equation}\label{3.33prim}
\nabla^\bot\theta(t)=-\nabla^\bot E*z(t)=-K(z(t)),\mbox{ a.e. }t\in(0,T),\end{equation}where $E(x)\equiv\frac1{2\pi}\ln|x|.$  
It follows, therefore, by \eqref{3.28} that
\begin{equation}
\label{3.28aaaa}
\barr{r}
K_\vp(z(t))=\nabla^\bot\theta_\vp(t)\to
-\nabla^\bot E*z(t)=-K(z(t)) \mbox{ weakly in } L^2\vsp\mbox{ for a.e. $t\in(0,T)$.}\earr\end{equation}
Now, by \eqref{e3.23a}, it follows by the lower semicontinuity of the $L^2$-norm that
$$|K(z(t))|_2\le C,\ \ \mbox{ a.e. }t\in(0,T),$$ and so $|K(z)|_2\in L^\9(0,T)$.  

For $0<\vp'<\vp\le1$ by the resolvent equation for $(\vp I-\Delta)\1$ and \eqref{3.28aaaa}, we have for a.e. $t\in(0,T)$ and $h\in L^2$, $|h|_2\le1$,
$$(h,K_\vp(z(t)))_2 
=(h,\nabla^\bot\Phi_{\vp'}(z(t)))_2
+\dd\frac{(\vp'{-}\vp)}\vp\,(h,\vp(\vp I{-}\Delta)\1\nabla^\bot\Phi_{\vp'}(z(t)))_2.$$
Hence, 
$$\barr{ll}
|(h,K_\vp(z(t)))_2|
&\le\dd\limsup_{\vp'\to0}
|(h,K_{\vp'}(z(t)))_2|\vsp
&+\dd\limsup_{\vp'\to0}
|(\vp(\vp I-\Delta)\1h,K_{\vp'}(z(t)))_2|\vsp 
&=|(h,K(z(t)))_2|+|(\vp(\vp I-\Delta)\1h,K(z(t)))_2|\vsp 
&\le2|K(z(t))|_2\ \mbox{ a.e. } t\in(0,T).\earr$$Therefore,
\begin{equation}\label{e3.17}
|K_\vp(z(t))|_2\le2|K(z(t))|_2,\ \mbox{ a.e. } t\in(0,T),\ \ff\vp>0.\end{equation}
We come back to \eqref{e3.15} and, taking into account \eqref{e3.4a} and \eqref{2.6secund}, we obtain that 
\begin{equation}\label{e3.25a}
\hspace*{-2mm}\barr{ll}
	2|K(z(t))u_1(t)|_2
	|K_\vp(z(t))|_2
	\le|K(z(t))u_1(t)|^2_2+|K_\vp(z(t))|^2_2\vsp
	\qquad\quad\le C|K(z(t))|_2|\nabla K(z(t))|_2|u_1(t)|^2_4+|K_\vp(z(t))|^2_2\vsp
	\qquad\quad\le C|K(z(t))|_2|z(t)|_2|u_1(t)|^2_4+|K_\vp(z(t))|^2_2\vsp
	\qquad\quad\le\dd\frac\nu4|z(t)|^2_2{+}\frac C\nu|K(z(t))|^2_2|u_1(t)|^4_4{+}|K_\vp(z(t))|^2_2,\, \mbox{a.e. }t\in(0,T).	
	\earr
\end{equation}
Then, substituting \eqref{e3.15aa}, \eqref{e3.25a}, \eqref{e3.17} in \eqref{e3.15}, we get 
\begin{equation}\label{e3.19}
\barr{r}
h'_\vp(t){+}|z(t)|^2_2\le C(h_\vp(t){+}(1{+}|u_1(t)|^4_4{+}|u_2(t)|^4_{\frac43})|K(z(t))|^2_2),\vsp \mbox{a.e. }t\in(0,T),\earr\end{equation}
Then, recalling that $|K(z)|_2\in L^\9(0,T)$ and that $z,u_1,u_2\in L^4(0,T;L^{\frac43}\cap L^4)$, it follows that the right-hand side of \eqref{e3.19} is in $L^1(0,T)$. Then, integrating over $(0,t)$ and taking into account \eqref{e3.14}, it follows via Gronwall's lemma that
$$h_\vp(t)\le C\int^t_0\(1+|u_1(s)|^4_4+|u_2(s)|^4_{\frac43}\)|K(z(s))|^2_2ds,\ \ff t\in(0,T),$$and, therefore, once again by  \eqref{e3.12} we have
$$|K_\vp(z(t))|^2_2\le C\int^t_0
\(1+|u_1(s)|^4_4+|u_2(s)|^4_{\frac43}\)|K(z(s))|^2_2ds, \mbox{ a.e. } t\in(0,T).$$ 
Then, by \eqref{3.28aaaa} and the weak lower semicontinuity of the $L^2$-norm, it follows that
$$|K(z(t))|^2_2\le C\int^t_0\(1+|u_1(s)|^4_4+|u_2(s)|^4_{\frac43}\)|K(z(s))|^2_2ds,\mbox{ a.e. } t\in(0,T).$$
The latter implies via Gronwall's lemma that $K(z(t))=0,$ for a.e. $t\in(0,T)$. Hence, by \eqref{3.33prim} and \eqref{e3.30prim}, 
we have therefore $z\equiv0$, a.e. on $[0,T]$, as claimed.\hfill$\Box$\bk

Theorem \ref{t2.2}  implies  the uniqueness of distributional solutions $y$ to \eqref{e1.7} in the sense of \eqref{e2.2}--\eqref{e2.3}  satisfying \eqref{e2.6}--\eqref{e2.7b} with initial data $u_0\in\calm(\rr^2)$, $u_0\ge0$. Namely, we have

\begin{corollary}\label{c3.1} Let $u_0\!\in\!\calp^a$  and let $u_1,u_2\in L^4(0,T;L^4)$ be two nonnegative distributional solutions to \eqref{e1.7} in the sense of \eqref{e2.2}, \eqref{e2.3}, such that $u_1-u_2\in L^\9(0,T;H\1)$. 
	Then $u_1\equiv u_2$.
	\end{corollary}

\n{\bf Proof.} We note first that, if $u$ is such a solution to \eqref{e1.7}, then we have
\begin{equation}\label{e3.23}
\int_{\rr^2}u(t,x)dx=\int_{\rr^2}u_0(dx)=1,\ \mbox{ a.e. }t\in(0,T).\end{equation}
  Indeed, if $u^\vp=u*\rho_\vp$, where $\rho_\vp=\frac1{\vp^3}\,\rho\(\frac t\vp,\frac x\vp\)$ is a mollifier  in $\rr^3$, then we have
  $$u^\vp_t-\nu\D u^\vp+\divv((K(u)u)*\rho_\vp)=0\mbox{ on }(0,T)\times\rr^2,$$and, integrating over  $\rr^2$, we get
  $$\int_{\rr^2}u^\vp(t,x)dx=\int_{\rr^2}u^\vp_0dx,\ \ff t\in(0,T),\ \ff\vp>0,$$which for $\vp\to0$ yields \eqref{e3.23}, as claimed.   
  
  In particular, $u_1,u_2\in L^\9(0,T;L^1)$, hence by interpolation $u_1,u_2\in L^4(0,T;L^4\cap L^{\frac43})$. 

Then, it follows from Lemma 2.3 in \cite{26secund} that there is a $dt\otimes dx$ version $\wt u$ of $u$ such that, for $\wt u(t,dx):=\wt u(t,x)dx,$ $t>0$, and $\wt u(0,dx):=u_0(dx)$ we have that $t\mapsto\int_{\rrd}\vf(x)\wt u(t,dx)$ is   continuous on $[0,T]$, for every $\vf\in C_b(\rr^2)$.

If $\wt u_1,\wt u_2$ are two such $dt\otimes dx$ versions of $u_1,u_2$, respectively, we have
$$\barr{l}
\dd\lim_{t\to0}\ {\rm ess}\sup_{\hspace*{-5mm}s\in(0,t)}|(u_1(s)-u_2(s),\vf)_2|\\\qquad
=\dd\lim_{t\to0}\ {\rm ess}\sup_{\hspace*{-5mm}s\in(0,t)}\left|\int_{\rr^2}(\wt u_1(s,x)-\wt u_2(s,x))\vf(x)dx\right|=0,\ \
\ff\vf\in C^\9_0(\rr^2).\earr$$
Then, by Theorem \ref{t2.2} it follows that 
$u_1\equiv u_2\, dt\otimes dx$ a.e., as claimed.\hfill$\Box$

\begin{remark}\label{r3.2}\rm Following the previous proof, one can get the uniqueness in Theorem \ref{t2.2} in the class of solutions $u\in L^{q_1}(0,T;L^{p_1})\cap L^{q_2}(0,T;L^{p_2})$, where $p_1\in(2,\9),$ $p_2\in(1,2)$ and $\frac1{q_1}+\frac1{p_1}\le1,$ $i=1,2$. This class of solutions seems to be more appropriate  if one takes into account   estimate \eqref{e2.4}.
\end{remark}

\section{Existence and weak uniqueness\\ for the McKean--Vlasov equation \eqref{e1.8}}\label{s4}
\setcounter{equation}{0}

It is well known (see \cite{2}, \cite{3}) that the existence of a distributional narrowly continuous solution $u$ for a \FP\ equation with Nemytskii drift terms implies the existence of a probabilistically weak solution $X(t)$ to the corresponding McKean--Vlasov equation such that 
\begin{equation}\label{e4.1}
	u(t,x)=\call_{X(t)}(x),\ \ff\,t>0\mbox{ and }u_0(dx)=\mathbb{P}\circ X(0)\1(dx),
\end{equation}where $\call_{X(t)}$ is the density of $\mathbb{P}\circ X(t)\1$ w.r.t. Lebesgue measure. 
This result follows by the nonlinear superposition principle in \cite[Section 2]{2} (which in turn is derived from the linear superposition principle in \cite{13}).    
Applying this to equation \eqref{e1.8} and, respectively, to the vorticity equation \eqref{e1.7}, which as seen earlier can be viewed as a \FP\ equation with the drift $K(u)$ satisfying $uK(u)\in L^1$,   the following existence result for equation \eqref{e1.8} follows by Theorem  \ref{t2.2bis} and \cite[Section 2]{2}. 

 \begin{theorem}\label{t4.1} Let $u_0\in\calp$.   Then,   there is a probabilistically weak solution $X$ to \eqref{e1.8} such that \eqref{e4.1} holds, where $u$ is the mild solution to equation \eqref{e1.7}, from Theorem {\rm\ref{t2.1}}. {\rm(See also Remark \ref{r2.1prim}.)}	
 \end{theorem} 

We recall that the process $X=X(t)$ is called a {\it probabilistically weak solution} to \eqref{e1.8} if there is a 2-dimensional $(\calf_t)$-Brownian  motion $W(t),$ $\ge0$, on a stochastic basis $(\ooo,\calf,(\calf_t)_{t\ge0},\mathbb{P})$ such that $X:[0,\9)\times\ooo\to\rr^2$ is progressively measurable, $\mathbb{P}$-a.s. continuous in $t$  and satisfies \eqref{e1.8}, i.e.,
\begin{equation}\label{e4.2}
dX(t)=K(u(t,\cdot))(X(t))dt+\sqrt{2\nu}\ dW(t),\ t\ge0,\end{equation}with one dimensional time  marginal laws  $\call_{X(t)}=\mathbb{P}\circ X(t)\1=u(t),\ t\ge0.$

The process $X(t)$ is called the {\it probabilistic representation} of the solution $u$ to the vorticity equation \eqref{e1.7}.  

In particular, we have by the Biot--Savart formula \eqref{e1.6} the probabilistic representation 
\begin{equation}\label{e4.3}
y(t)=K(\call_{X(t)}),\ \ff t\ge0,\end{equation}of the solution $y$ to the Navier--Stokes equation \eqref{e1.1}. 

We shall discuss now the weak uniqueness of probabilistically weak solutions $X$ to the McKean--Vlasov equation \eqref{e1.8}. To this purpose, we shall prove first the  {\it linearized uniqueness} for equation \eqref{e1.7}. Namely, we have

\begin{theorem}\label{t4.2} Let $u\in L^4(0,T;L^{\frac43})$ and let
\begin{equation}\label{e4.4}
\barr{c}
u_1,u_2\in L^4(0,T;L^4\cap L^{\frac43}),\ u_1-u_2\in L^\9(0,T;H\1),\vsp u_1(t)-u_2(t)\in L^1\mbox{ for a.e. }t\in(0,T),\earr \end{equation}such that 
\begin{equation}\label{e4.4a}
	\lim_{t\downarrow 0}\,
	{\rm ess}\sup_{\hspace*{-6mm}{s\in(0,T)}}
\int_{\rr^2}(u_1(s,x)-u_2(s,x))\vf(x)dx=0,\ \ff\vf\in C^\9_0(\rr^2)
\end{equation}
be two solutions to the equation
\begin{equation}\label{e4.5}
\barr{r}\dd
\int^\9_0\int_{\rr^2}(\vf_t+\nu\D\vf +K(u)\cdot\nabla\vf)v\,dxdt
+\int_{\rr^2}\vf(0,x)u_0(dx)=0,\\ \ff\vf\in C^2_0([0,T];\rr^2),\earr\end{equation}where $u_0\in\calm(\rr^2)$. Then, $u_1\equiv u_2.$	
\end{theorem}

\n{\bf Proof.}  The proof is the same as that of Theorem \ref{t2.2}, but with some simplifications. Namely, we set $z=u_1-u_2$ and get, by \eqref{e4.1}, 
\begin{equation}\label{e4.6}
\barr{l}
z_t-\nu\D z+\divv(K(u)z)=0\mbox{ in }\cald'((0,T)\times\rr^2),\vsp 
z(0)=u_0.\earr\end{equation}If $z_\vp=\Phi_\vp(z)$, we obtain for $z_\vp$ the equation
$$(z_\vp)_t-\nu\D z_\vp+\Phi_\vp(\divv(K(u)z))=0\mbox{ in }(0,T)\times\rr^2$$and so, arguing  as in the proof of Theorem \ref{t2.2}, we obtain (see \eqref{e3.12}--\eqref{e3.13})
\begin{equation}\label{e4.7}
\barr{c}
|K_\vp(z)(t)|^2_2+\dd\int^t_0|z(s))|^2_2ds
\le C\dd\int^t_0|K(u(s))|^2_4|K_\vp(z)|^2_2ds\vsp 
\le C\dd\int^t_0|u(s)|^2_{\frac43}|K_\vp(z)|^2_2ds,\ \ff t\in[0,T],\earr\end{equation}from which, taking into account \eqref{e3.17}, we get for $\vp\to0$ 
$$|K(z(t))|^2_2\le C\int^t_0(1+|u(s)|^4_4+|u(s)|^4_{\frac43})|K(z(s))|^2_2ds,\mbox{ a.e. }t\in(0,T),$$which implies as above that $z\equiv0$.\hfill$\Box$\bk

Arguing as in the proof of Corollary \ref{c3.1}, it follows, by Theorem \ref{t4.2}, the following uniqueness result.

\begin{corollary}\label{c4.3} Let $u_0\in\calp^a$ and let $u_1,u_2\in L^4(0,T;L^4)$ be two nonnegative solutions to \eqref{e4.5} such that $u_1-u_2\in L^\9(0,T;H\1)$. Then, $u_1\equiv u_2$.
	\end{corollary}

Now, let us prove  weak uniqueness for the McKean--Vlasov equation \eqref{e1.8}. (For the definition of weak solutions we refer to Definition 3.1 (a) part (i) in \cite{27prim}.) 

\begin{theorem}\label{t4.3} Let $T>0$ and let $X(t)$, $\wt X(t),$ $t\ge0$, on stochastic bases $(\ooo,\calf,(\calf_t)_{t\ge0},\mathbb{P})$, $(\wt\ooo,\wt\calf,(\wt\calf_t)_{t\ge0},\wt{\mathbb{P}})$ respectively, be two probabilistically weak solutions to \eqref{e1.8} such that, for	 
	$$u(t,\cdot)=\call_{X(t)},\ \ \ \wt u(t,\cdot)=\call_{\wt X(t)},\ t>0,$$we have
\begin{equation}\label{e4.8}
u,\wt u\in L^4(0,T;L^4)\cap L^\9(0,T;H\1).	 
	\end{equation}Then $X$ and $\wt X$ have the same laws, that is,
\begin{equation}\label{e4.9}
	\mathbb{P}\circ X\1=\wt{\mathbb{P}}\circ\wt X\1.	\end{equation}
\end{theorem}

\n{\bf Proof.}  By It\^o's formula, both $u$ and $\wt u$ satisfy \eqref{e2.3} and by our definition of weak solution also \eqref{e2.2}. Hence, by Corollary \ref{c3.1}, $u\equiv\wt u$. Furthermore, again by It\^o's formula, both $\mathbb{P}\circ X\1$ and $\wt{\mathbb{P}}\circ\wt X^{-1}$ satisfy the martingale problem with the initial condition $u_0$ for the linearized Kolmogorov operator
\begin{equation}\label{e4.11}
L_u:=\Delta+K(u)\cdot\nabla.\end{equation}Note that Theorem \ref{t4.2} above remains true if we replace the role of the starting time $0$ by any $s\ge0$. Therefore, the assertion follows by Lemma 2.12 in \cite{13} applied to the linear Kolmogorov operator in \eqref{e4.11} and the family $\calr_{[s,T]}$, $0\le s\le T$, of $\calr$-regular narrowly continuous solutions are defined as follows. $\calr_{[s,T]}$ is the set of all narrowly continuous solutions $(u(t))_{t\ge s}$ of \eqref{e2.2}, \eqref{e2.3} starting at $(s,\zeta)$ with $\zeta\in\calp(\rr^d)$ such that
$$u\in L^4(s,T;L^4)\cap L^\9(s,T;H\1).$$Obviously, this family fulfills conditions (2.9) and (2.14) in \cite{13}, as is required for Lemma 2.12 in \cite{13}.\hfill$\Box$\bk
 
Now, let us turn to the probabilistically strong solutions to \eqref{e1.8}. 

\begin{theorem}\label{t4.5} Let $u_0\in\calp^a\cap L^4$.  Then, the solution to \eqref{e1.8} from Theorem {\rm\ref{t4.1}} is, in fact, a  probabilistically strong solution, i.e., is a functional of the Brownian motion $W(t)$, $t\ge0$. Furthermore, pathwise uniqueness holds in the class of all probabilistically weak  solutions to \eqref{e1.8} with the same Brownian motion, having path laws with one dimensional time marginal law densities in $L^{\frac43}(0,T;L^{\frac43})$.   
	\end{theorem}

\n{\bf Proof.} Let $u$ be as in Theorem \ref{t4.1}, with initial $u_0$ and let $(X,W)$ be the corresponding weak solution to \eqref{e1.8}. Then, fixing $u$ in \eqref{e1.8}, we are in the case of a usual SDE and may apply the results in Sections 1.3 and 2.1.1 in \cite{21prim}. To this end, we aim to prove 
\begin{equation}\label{4.12}
K(u)\in L^\9(0,T;W^{1,4}).\end{equation}
Since $u_0\in\calp^a\cap L^4\subset L^{\frac43},$ by Theorem \ref{t2.2bis}~(iii) we know that $u\in L^\9(0,T;L^{\frac43})$, hence by \eqref{e1.6a} and \eqref{2.12prim} we have $K(u)\in L^\9(0,T;L^4)$. Furthermore, by \eqref{2.6secund} and \eqref{2.12prim}, $\nabla K(u)\in L^\9(0,T;L^4)$ and thus \eqref{4.12} is proved.  Hence, the assertion follows by Theorems 1.3.1 and 2.1.3 in \cite{21prim} and Lemmas A.2 and A.3 in \cite{16prim}.\hfill$\Box$\bk

Finally, we prove that the path laws of the probabilistically weak solution to \eqref{e1.8} from Theorem \ref{t4.1} form a Markov process. To this end, we first note that clearly both Theorems \ref{t2.1} and \ref{t2.2} hold if we consider \eqref{e1.7} on $[s,\9)\times\rr^2$ for any $s\ge0$. Then, renaming the initial condition in \eqref{e2.3} by $\zeta\in\calp$, due to Theorem \ref{t2.1} and Remark \ref{r2.1prim} we have, for each $(s,\zeta)\in [0,\9)\times\calp,$ a solution to \eqref{e2.2}, \eqref{e2.3} with the initial condition $\zeta$ at time $s$, which according to Lemma 2.3 in \cite{26secund} has a narrowly continuous version  on $[0,\9)$. 
Let us denote this narrowly continuous solution by $\mu^{s,\zeta}=(\mu^{s,\zeta}_t)_{t\ge s}$. Below we identify a measure which is absolutely continuous w.r.t. Lebesgue measure $dx$ with its density. Furthermore, for $(s,\zeta)\in[0,\9)\times\calp$, we denote the corresponding probabilistically weak solution by $X(t,s,\zeta)_{t\ge s}$, defined on a stochastic basis $(\ooo^{s,\zeta},\calf^{s,\zeta},(\calf^{s,\zeta})_{t\ge s},\mathbb{P}^{s,\zeta})$ and define
$$\mathbb{P}_{s,\zeta}:=\mathbb{P}^{s,\zeta}\circ X(\cdot,s,\zeta)\1.$$Then, $\mathbb{P}_{s,\zeta}$ is a probability measure on $\ooo_s:=C([s,\9);\rr^d)$, i.e., the set of all continuous paths in $\rr^d$ starting at time $s$  equipped with the topology of locally uniform convergence and corresponding Borel $\sigma$-algebra $\calb(\ooo_s)$. 

Define, for $\tau\ge s$,
$$\pi^s_\tau:\ooo_s\to\rr^d,\ \ \pi^s_\tau(w):=w(\tau),\ w\in\ooo_s,$$and, for $r\ge s$,
$$\calf_{s,r}:=\sigma(\pi^s_\tau:s\le\tau\le r).$$

\begin{theorem}
	\label{t4.6} The family $\mathbb{P}_{s,\zeta},\ (s,\zeta)\in[0,\9)\times\calp$, forms a nonlinear Markov process in the sense of Definition {\rm\ref{d4.7}} below with $\calp_0:=\calp$.
	\end{theorem}

The following is a moderately concretized version of the one by McKean from \cite{24prim}. 

\begin{definition}\label{d4.7} \rm 	Let $\calp _0 \subseteq \calp $. A \textit{nonlinear Markov process} is a family\break $(\mathbb{P}_{s,\zeta})_{(s,\zeta)\in \rr_+\times \calp _0}$ of probability measures $\mathbb{P}_{s,\zeta}$ on $\calb (\Omega_s)$ such that
	\begin{enumerate}
		\item[(i)] The marginals $\mathbb{P}_{s,\zeta}\circ(\pi^s_t)^{-1} =: \mu^{s,\zeta}_t$ belong to $\mathcal{P}_0$ for all $0\leq s \leq r \leq t$ and $\zeta \in \calp _0$.
		\item[(ii)] The \textit{nonlinear Markov property} holds, i.e. for all $0\leq s \leq r \leq t$, $\zeta \in \calp _0$
		\begin{equation}\label{Markov-prop}\tag{MP}
			\barr{r}
			\mathbb{P}_{s,\zeta}(\pi^{s}_{t} \in A|\mathcal{F}_{s,r})(\cdot) = p_{(s,\zeta),(r,\pi^s_r(\cdot))}(\pi^r_t\in A) \quad \mathbb{P}_{s,\zeta}-\text{a.s.}\vsp \text{ for all }A \in \calb (\rr^d),\earr
		\end{equation}
		where $p_{(s,\zeta),(r,y)}, y \in \rr^d$, is a regular conditional probability kernel from $\rr^d$ to $\calb (\Omega_r)$ of $\mathbb{P}_{r,\mu^{s,\zeta}_r}[\,\,\cdot\, \,| \pi^r_r=y]$, $y \in \rr^d$ (i.e., in particular,  $p_{(s,\zeta),(r,y)} \in \calp (\Omega_r)$ and $p_{(s,\zeta),(r,y)}(\pi^r_r = y) = 1$).
	\end{enumerate}
\end{definition}

The term \emph{nonlinear} Markov property originates from the fact that in the situation of Definition \ref{d4.7} the map $\calp _0 \ni \zeta \mapsto \mu^{s,\zeta}_t$ is, in general, not convex.
	
 \begin{remark}\label{r4.8} \rm The one-dimensional time marginals $\mu^{s,\zeta}_t = \mathbb{P}_{s,\zeta}\circ (\pi^s_t)^{-1}$ of a nonlinear Markov process satisfy the {\it flow property}, i.e.,
 	\begin{equation}\label{4.13}
 		\mu^{s,\zeta}_t = \mu^{r,\mu^{s,\zeta}_r}_t,\quad \forall\, 0\leq s \leq r \leq t, \zeta \in \calp_0.
 	\end{equation}
 
\end{remark}

\n{\bf Proof of Theorem \ref{t4.6}.} For $(s,\zeta)\in[0,\9)\times\calp$, consider the narrowly continuous solution $\mu^{s,\zeta}=(\mu^{s,\zeta}_t)_{t\ge s}$, to \eqref{e2.2}, \eqref{e2.3} introduced in front of the formulation of Theorem \ref{t4.6} and define
$$\textgoth{P}_0:=\{u_0\in\calp^a:u_0\in L^4\}.$$ 
Then, by Theorem \ref{t2.2bis}~(iii) we have for every $(s,u_0)\in [0,\9)\times\textgoth{P}_0$,
$$\mu^{s,u_0}\in L^\9(s,T;L^4)\subset L^4(s,T;L^4)\cap L^\9(s,T;H\1),$$and \eqref{e2.5a}, \eqref{2.7} imply that $\mu^{s,u_0}_t\in\textgoth{P}_0$ for every $0\le s\le t$, $u_0\in\textgoth{P}_0$ and $\mu^{t,\zeta}\in\textgoth{P}_0$ for every $0\le s<t$, $\zeta\in\calp.$ Furthermore, Corollary \ref{c3.1} implies that $\mu^{s,u_0}$, $(s,u_0)\in[0,\9)\times \textgoth{P}_0,$ satisfy the flow property \eqref{4.13}. Then, by Corollary \ref{c4.3} and \cite[Lemma 3.4]{27prim} we see that Corollary 3.9 in \cite{27prim} (with $\calp_0:=\calp$) applies to prove the assertion.\hfill$\Box$



\appendix

\section*{Proof of \eqref{e3.23b}}
\setcounter{theorem}{0}

\begin{lemma}\label{A.1} Let $z\in L^{\frac43},\vp>0.$ Then
	\begin{equation}\label{eA.1}
	K_\vp(z)=-K(z)+\vp K(\Phi_\vp(z)).
	\end{equation}
\end{lemma}

\pf We first recall that $\vp\Phi_\vp$ is a contraction on every $L^p$, $p\in[1,\9]$, and that
$$\Phi_\vp(z)=(g_\vp*z),$$where
\begin{equation}\label{eA.2}
g_\vp(x):=\int^\9_0 e^{-\vp t}\ \frac1{4\pi t}\ e^{-\frac1{4t}\,|x|^2}dt,\ x\in\rr^2. 
\end{equation}(See, e.g., \cite[p.~132, formula (26)]{24p}.) Then, $\vp g_\vp\in\calp^a$ and an elementary  com\-pu\-ta\-tion yields
\begin{equation}\label{eA.3}
\nabla^\bot g_\vp(x)=-k(x)\int^\9_0 e^{-\vp|x|^2t}\,\gamma(dt),\ x\in\rr^2\setminus\{0\},
\end{equation}where $\gamma$ is a probability measure on $[0,\9)$ with density
$$t\mapsto\frac1{4t^2}\ e^{-\frac1{4t}},\ t\in[0,\9),$$ and $k$ is as in \eqref{e2.6'}. Hence, if $B_1$ denotes the unit ball in $\rr^2$ with centre zero, we have
\begin{equation}\label{eA.4}
\dd\sup_{\vp>0}\(\one_{B_1}|\nabla^\bot g_\vp|\)\le\one_{B_1}|k|\in L^1\mbox{\ \ and\ \ }
\dd\sup_{\vp>0}\(\one_{B^c_1}|\nabla^\bot g_\vp|\)\le\one_{B^c_1}|k|\in L^\9
\end{equation}and
\begin{equation}\label{eA.5}
|\nabla^\bot g_\vp(x)|\nearrow|k(x)|=\frac1{2\pi|x|},\ \ \ff x\in\rr^2\setminus\{0\}.
\end{equation}We first prove \eqref{eA.1} for $\vf\in C^\9_0=:\cald.$ Let $\cald'$ denote its dual and $\cals$ the space of all Schwartz test functions. Then, we have by the resolvent equation of $\Phi_{\vp'}$, $\vp'>0$, that, for all $\vp'\in(0,\vp)$,
\begin{equation}\label{eA.6}
\Phi_\vp(\vf)=\Phi_{\vp'}(\vf)+(\vp'-\vp)\Phi_{\vp'}(\Phi_\vp(\vf)).
\end{equation}
By \eqref{eA.4}, \eqref{eA.5} and Lebesgue's dominated convergence theorem, we have
$$\nabla^\bot\Phi_{\vp'}(\vf)
=(\nabla^\bot g_{\vp'})*\vf
\underset{\vp'\to0}{\longrightarrow}-k*\vf=-K(\vf)\mbox{ in }L^1,\mbox{ hence in }\cald',$$and, for all $\wt\vf\in\cald,$ since $\Phi_\vp(\vf)\in\cals$, hence $\Phi_\vp(\vf)*\wt\vf\in\cals$,
$$\barr{ll}
{}_{\raise-3pt\hbox{$_{\cald'}$}}\!\<\nabla^\bot\Phi_{\vp'}(\Phi_{\vp}(\vf)),
\wt\vf\>_{\cald}\!\!\!\! 
&=\dd\int_{\rr^2}\!\!\nabla^\bot g_{\vp'}\Phi_{\vp}(\vf)*\wt\vf\,dx 
\underset{\vp'\to0}{\longrightarrow}-\dd\int_{\rr^2}\!\!k\,\Phi_{\vp}(\vf)*\wt\vf\,dx\vsp&=-\dd\int_{\rr^2}\!\! k*\Phi_{\vp}(\vf)\wt\vf\,dx\vsp 
&=-\dd\int_{\rr^2}\!\! K(\Phi_{\vp}(\vf))\wt\vf\,dx.
\earr$$Therefore,
$$\nabla^\bot \Phi_{\vp'}(\Phi_{\vp}(\vf))\to- K(\Phi_{\vp}(\vf))\mbox{\ \ in }\cald',$$
and, consequently, applying $\nabla^\bot$ to \eqref{eA.6}, and passing to the limit in $\cald'$ with $\vp'\to0$, we obtain that \eqref{eA.1} holds for $\vf\in C^\9_0$. Now, we approximate $z$ in $L^{\frac43}$ by $\vf_n\in C^\9_0$, $n\in\nn$, and since then also $\Phi_\vp(\vf_n)\to\Phi_\vp(z)$ in $L^{\frac43}$ as $n\to\9$, using the generalized Young inequality and the fact that $k$, $\nabla^\bot g_\vp\in L^{2,\9}$, we can pass to the limit with $n\to\9$ in $L^4$ to obtain \eqref{eA.1}, for $z\in L^{\frac43}$.\hfill$\Box$

 
\bk\n{\bf Acknowledgement.} M. R\"ockner was supported by the DFG through SFB 1283/2 2021-317210226, V. Barbu was supported by the  CNCS--UEFISCDI project  PN-III-P4-PCE-2021-0006, within PNCDI III and D. Zhang by\break the NSFC (No. 12271352, 12322108) and the Shanghai Rising Star Program\break 21QA1404500.

\end{document}